# Fuzzy Statistical Limits

Mark Burgin[a] and Oktay Duman[b]

[a]*Department of Mathematics, University of California, Los Angeles, California 90095-1555, USA*

[b]*TOBB Economics and Technology University, Faculty of Arts and Sciences, Department of Mathematics, Söğütözü 06530, Ankara, Turkey*

**Abstract**

Statistical limits are defined relaxing conditions on conventional convergence. The main idea of the statistical convergence of a sequence *l* is that the majority of elements from *l* converge and we do not care what is going on with other elements. At the same time, it is known that sequences that come from real life sources, such as measurement and computation, do not allow, in a general case, to test whether they converge or statistically converge in the strict mathematical sense. To overcome these limitations, fuzzy convergence was introduced earlier in the context of neoclassical analysis and fuzzy statistical convergence is introduced and studied in this paper. We find relations between fuzzy statistical convergence of a sequence and fuzzy statistical convergence of its subsequences (Theorem 2.1), as well as between fuzzy statistical convergence of a sequence and conventional convergence of its subsequences (Theorem 2.2). It is demonstrated what operations with fuzzy statistical limits are induced by operations on sequences (Theorem 2.3) and how fuzzy statistical limits of different sequences influence one another (Theorem 2.4). In Section 3, relations between fuzzy statistical convergence and fuzzy convergence of statistical characteristics, such as the mean (average) and standard deviation, are studied (Theorems 3.1 and 3.2).

*Keywords*: statistical convergence, fuzzy sets, fuzzy limits, statistics, mean, standard deviation, fuzzy convergence

## 1. Introduction

Statistical limits are defined relaxing conditions on conventional convergence. Conventional convergence in analysis demands almost all elements of the sequence to satisfy the convergence condition. Namely, almost all elements of the sequence have to belong to arbitrarily small neighborhood of the limit. The main idea of statistical convergence is to relax this condition and to demand validity of the convergence condition only for a majority of elements. As always, statistics cares only about big quantities and majority is a surrogate of the concept "almost all" in pure mathematics. The reason is that statistics works with finite populations and samples, while pure mathematics is mostly interested in infinite sets.

However, the idea of statistical convergence, which emerged in the first edition (published in Warsaw in 1935) of the monograph of Zygmund [19], stemmed not from statistics but from problems of series summation. The concept of statistical convergence was formalized by Steinhaus [18] and Fast [10] and later reintroduced by Schoenberg [16]. Since that time, statistical convergence has become an area of active research. Researchers studied properties of statistical convergence and applied this concept in various areas: measure theory [13], trigonometric series [19], approximation theory [9], locally convex spaces [12], summability theory and the limit points of sequences [7], [15], densities of natural numbers [14], in the study of subsets of the Stone-Čhech compactification of the set of natural numbers [6], and Banach spaces [8].

However, as a rule, neither limits nor statistical limits can be calculated or measured with absolute precision. To reflect this imprecision and to model it by mathematical structures, several approaches in mathematics have been developed: fuzzy set theory, fuzzy logic, interval analysis, set valued analysis, etc. One of these approaches is the neoclassical analysis (cf., for example, [2, 3]). In it, ordinary structures of analysis, that is, functions, sequences, series, and operators, are studied by means of fuzzy concepts: fuzzy limits, fuzzy continuity, and fuzzy derivatives. For example, continuous functions, which are studied in the classical analysis, become a part of the set of the fuzzy continuous functions studied in neoclassical analysis. Neoclassical analysis extends methods of classical calculus to reflect uncertainties that arise in computations and measurements.

The aim of the present paper is to extend and study the concept of statistical convergence utilizing a fuzzy logic approach and principles of the *neoclassical analysis*, which is a new branch of fuzzy mathematics and extends possibilities provided by the classical analysis [2, 3]. Ideas of fuzzy logic and fuzzy set theory have been used not only in many applications, such as, in bifurcation of non-linear dynamical systems, in the control of chaos, in the computer programming, in the quantum physics, but also in various branches of mathematics, such as, theory of metric and topological spaces, studies of convergence of sequences and functions, in the theory of linear systems, etc.

In the second section of this paper, going after introduction, we introduce a new type of statistical convergence, the concept of fuzzy statistical convergence, and give a useful characterization of this type of convergence. In the third section, we consider relations between fuzzy statistical convergence and fuzzy convergence of statistical characteristics such as the mean (average) and standard deviation.

For simplicity, we consider here only sequences of real numbers. However, in a similar way, it is possible to define statistical fuzzy convergence for sequences of complex numbers and obtain similar properties.

In what follows, *N* denotes the set of all natural numbers and *R* denotes the set of all natural numbers.

## 2. Fuzzy statistical convergence

Here we extend the concept of statistical convergence to the concept of fuzzy statistical convergence, which, as we have discussed, is more realistic for real life applications.

For convenience, throughout the paper, $r$ denotes a non-negative real number and $l = \{a_i\,;\,i = 1, 2, 3, \ldots\}$ represents a sequence of real numbers.

Consider a subset $K$ of the set $N$. Let $K_n = \{k \in K;\, k \leq n\}$.

The (*asymptotic*) *density* $d(K)$ of a set $K$ is equal to

$$\lim_{n \to \infty} (1/n) |K_n|$$

whenever this limit exists (cf., for example, [11], [14]).

Changing the conventional limit in this formula to the *r*-limit [3], we obtain fuzzy asymptotic density.

**Definition 2.1.** The *asymptotic r-density* $d_r(K)$ of the set $K$ is equal to

$$r\text{-}\lim_{n \to \infty} (1/n) |K_n|,$$

whenever the fuzzy limit exists; here $|B|$ denotes the cardinality of the set $B$.

As a rule, *r*-limit of a sequence is not uniquely defined, i.e., one sequence can have many *r*-limits. Consequently, several numbers can be equal to the asymptotic *r*-density $d_r(K)$ of the same set $K$.

**Example 2.1.** Let us take $K = \{ 2i;\ i = 1, 2, 3, \ldots \}$. Then $½ = d_{½}(K)$, $¼ = d_{½}(K)$, and $3/5 = d_{½}(K)$.

However, in some cases, the asymptotic *r*-density $d_r(K)$ is unique.

**Lemma 2.1.** The asymptotic 0-density of any subset $K$ of $N$ (if it exists) coincides with the asymptotic density of $K$.

This shows that asymptotic *r*-density is a natural extension of the concept of asymptotic density.

**Lemma 2.2.** If $x = d_r(K)$, then $x = d_q(K)$ for any $q > r$.

**Corollary 2.1.** If the asymptotic density $d(K)$ of a subset $K$ of $N$ is equal to $x$, then $x = d_r(K)$ for any $r > 0$.

However, the following example shows that the converse of Corollary 2.1 is not always true. Moreover, asymptotic $r$-density can exist in such cases when asymptotic density does not exist.

**Example 2.2.** Let $K$ be the set of all even positive integers with an even number of digits to base ten. Then it is known from [14; pp: 248-249] that the density $d(K)$ does not exists. However, Definition 2.1 directly implies $d_{½}(K) = 0$.

Furthermore, as $0 \leq |K_n| \leq n$, Definition 2.1 yields the following results.

**Lemma 2.3.** If $k = r\text{-lim}_{n \to \infty} (1/n) |K_n|$ for a set $K \subseteq N$, then $-r \leq k \leq 1 + r$.

**Lemma 2.4.** $0 = 1\text{-lim}_{n \to \infty} (1/n) |K_n|$ and $1 = 1\text{-lim}_{n \to \infty} (1/n) |K_n|$ for any set $K \subseteq N$.

It means that for $r = 1$, the concept of asymptotic $r$-density, in some sense, degenerates as any number between 0 and 1 becomes the asymptotic $r$-density of any set $K \subseteq N$.

Let $K$ and $H$ be subsets of $N$, then properties of fuzzy limits imply the following result.

**Lemma 2.5.** If $k = d_r(K)$ and $h = d_q(H)$, then $d_{r+q}(K \cup H)$ exists and
$$d_{r+q}(K \cup H) \leq k + h$$

**Corollary 2.2.** If $k = d(K)$ and $h = d(H)$, then $d(K \cup H)$ exists and
$$d(K \cup H) \leq k + h$$

**Corollary 2.3.** If $k = d_r(K)$ and $h = d(H)$, then $d_r(K \cup H)$ exists and
$$d_r(K \cup H) \leq k + h$$

As in the case of asymptotic density [5], there are connections between asymptotic *r*-density and fuzzy Cesáro summability.

As usual, the well-known Cesáro matrix $C = \{c_{nk} ; n,k = 1,2,3,...\}$ is defined by its coefficients

$$c_{nk} = \begin{cases} 1/n, & \text{if } 1 \leq k \leq n \\ 0, & \text{otherwise.} \end{cases}$$

Using summability methods, we obtain the following result.

**Proposition 2.1.** If the asymptotic *r*-density $d_r(K)$ of a subset $K$ of $N$ exists, then we have

$$d_r(K) = r\text{-lim}_{n \to \infty} (1/n) \sum_{1 \leq k \leq n} \chi_K(k) = r\text{-lim}_{n \to \infty} (C\chi_K)_n$$

where $\chi_K$ denotes the characteristic function of the set $K$, and $C\chi_K$ denotes the sequence of Cesáro matrix transformation.

Let us consider a sequence $l = \{a_i ; i = 1, 2, 3, \dots \}$ of real numbers, real number $a$, and the set $L_\varepsilon(a) = \{i \in N; |a_i - a| \geq \varepsilon \}$.

**Definition 2.2.** The *asymptotic r-density*, or simply, *r-density* $d_r(l)$ of the sequence $l$ with respect to $a$ and $\varepsilon$ is equal to $d_r(L_\varepsilon(a))$.

**Lemma 2.6.** If $x = d_r(L_\varepsilon(a))$, then $x = d_q(L_\varepsilon(a))$ for any $q > r$.

This asymptotic density allows us to define statistical convergence.

**Definition 2.3.** A sequence $l$ *r-statistically converges* to a number $a$ if $0 = d_r(L_\varepsilon(a))$ for every $\varepsilon > 0$. The number (point) $a$ is called an *r-statistical limit* of $l$. We denote this limit by $a = r\text{-stat-lim } l$.

Then, Definition 2.3 implies the following result.

**Lemma 2.7.**

(a) $a = r\text{-stat-lim } l \Leftrightarrow \forall \, \varepsilon > 0$,

$$0 = r\text{-}\lim_{n \to \infty} (1/n) \, |\{i \in N; i \leq n \, ; |a_i - a| \geq \varepsilon\}|.$$

(b) $a = r\text{-stat-lim } l \Leftrightarrow \forall \, \varepsilon > 0$,

$$1 = r\text{-}\lim_{n \to \infty} (1/n) \, |\{i \in N; i \leq n \, ; |a_i - a| < \varepsilon\}|.$$

**Remark 2.1.** We know from [3] that if $a = \lim l$ (in the ordinary sense), then for any $r \geq 0$, we have $a = r\text{-lim } l$. In a similar way, using Definition 2.3, we can easily see that if $a = stat\text{-lim } l$, then we have $a = r\text{-stat-lim } l$ for any $r \geq 0$. However, its converse is not true. It is demonstrated in the following example of a sequence that is $r$-statistically convergent but not convergent and also not statistically convergent.

**Example 2.3.** Let us consider the sequence $l = \{a_i \, ; i = 1,2,3,\ldots\}$ whose terms are

$$a_i = \begin{cases} i, & \text{if } i = 2n \ (n = 1,2,3,\ldots) \\ 1, & \text{otherwise.} \end{cases}$$

Then, it is easy to see that the sequence $l$ is divergent in the ordinary sense. Even more, the sequence $l$ has no $r$-limits for any $r$ since it is unbounded from above (see Theorem 2.3 from [3]).

Likewise, the sequence $l$ is not statistically convergent.

At the same time, $1 = (\tfrac{1}{2})\text{-}stat\text{-lim } l$ since $0 = d_{\frac{1}{2}}(K)$ where $K = \{2n; n \in N\}$.

**Remark 2.2.** Example 2.3 also demonstrates a difference between the three types of convergence: conventional fuzzy and fuzzy statistical convergence. It is known that every convergent (in the ordinary sense) or *r*-convergent sequence is bounded. However, this situation is not valid for *r*-statistically convergent sequences.

Lemma 2.1 implies the following result.

**Lemma 2.8.** 0-statistical convergence coincides with statistical convergence.

This result shows that fuzzy statistical convergence is a natural extension of statistical convergence.

As in the case of statistical convergence [5], there are connections between fuzzy statistical convergence and fuzzy Cesáro summability. Namely, Proposition 2.1 and Definition 2.3 imply the following result.

**Proposition 2.2.** $a = r\text{-stat-lim } l$ if and only if $0 = r\text{-lim}_{n\to\infty}(C\chi_{L\varepsilon(a)})_n$

Using the definitions of statistical convergence and *r*-statistical convergence, it is possible to prove the following result.

**Lemma 2.9.** If $a = r\text{-stat-lim } l$, then $a = q\text{-stat-lim } l$ for any $q > r$.

**Corollary 2.3.** If $a = stat\text{-lim } l$, then $a = r\text{-stat-lim } l$ for any $r \geq 0$.

**Proposition 2.3.** Any sequence *l* is 1-statistically convergent.

Indeed, by Lemma 2.4, $0 = 1\text{-lim}_{n\to\infty}(1/n)|K_n|$ and $1 = 1\text{-lim}_{n\to\infty}(1/n)|K_n|$ for any set $K \subseteq N$. Consequently, $0 = d_1(L_\varepsilon(0))$ for every $\varepsilon > 0$ and any sequence *l* of real numbers as the number of elements in the set $L_{n\varepsilon}(0) = \{k \in L_\varepsilon(0); k \leq n\}$ is always less than *n*. Thus, the sequence *l* 1-statistically converges to zero.

The result of Proposition 2.3 shows that *r*-statistical convergence is interesting only for small *r*.

Let us take a sequence $l = \{a_i\, ;\, i = 1, 2, 3, \ldots\}$ and its sequence *h*.

**Definition 2.4.** The *asymptotic density*, or simply, *density* $d_l(h)$ of the subsequence *h* of the sequence *l* is equal to $d(\{i \in N;\, a_i \in h\})$.

**Lemma 2.10.** If *h* is a subsequence of the sequence *l* and *k* is a subsequence of the sequence *h*, then $d_l(k) \leq d_l(h)$.

Note that it is possible that $d_h(k) > d_l(h)$.

**Definition 2.5.** The *asymptotic r-density*, or simply, *r-density* $d_{r,l}(h)$ of a subsequence *h* of the sequence *l* is equal to $d_r(\{i \in N;\, a_i \in h\})$.

**Lemma 2.11.** If *h* is a subsequence of the sequence *l* and *k* is a subsequence of the sequence *h*, then $d_{r,l}(k) \leq d_{r,l}(h)$ for any $r > 0$.

Note that it is possible that $d_{r,h}(k) > d_{r,l}(h)$.

**Definition 2.6** [4]**.** A subsequence *h* of a sequence *l* is called *statistically dense* in *l* if $d_l(h) = 1$.

It is known that a subsequence of a fuzzy convergent sequence is fuzzy convergent [3]. However, for statistical convergence this is not true. Indeed, the sequence $h = \{i\, ;\, i = 1, 2, 3, \ldots\}$ is a subsequence of the fuzzy statistically convergent sequence *l* from Example 2.3. At the same time, *h* is statistically fuzzy divergent.

However, if we consider statistically dense subsequences of statistically fuzzy convergent sequences, it is possible to prove the following result.

**Theorem 2.1.** A sequence is *r*-statistically convergent if and only if any its statistically dense subsequence is *r*-statistically convergent.

<u>Proof</u>. Since the sufficiency is obvious, it is enough to prove the necessity. So, assume that a sequence $l = \{a_i\,;\,i = 1, 2, 3,…\}$ is *r*-statistically convergent, say $a = r$-*stat*-lim $l$, and that $h$ is any statistically dense subsequence of $l$. We will show that $a = r$-*stat*-lim $h$. By hypothesis, we may write for every $\varepsilon > 0$, there exists an increasing index subset $K$ of $N$ with $d_r(K)=1$ such that $|a_i - a| < \varepsilon$ for sufficiently large $i \in K$ (see Lemma 2.6.b). Also, since $h$ is any statistically dense subsequence of $l$, there exists an increasing index subset $M$ with $d(M)=1$ such that $h = \{a_i\,;\,i \in M\}$. Now we can see that $d_r(K \cap M) = 1$, and so we have $|a_i - a| < \varepsilon$ for sufficiently large $i \in K \cap M$. This implies that $a = r$-*stat*-lim $h$ and completes the proof.

An *r*-statistically convergent sequence contains not only dense *r*-statistically convergent subsequences, but also dense convergent subsequences.

The following theorem gives us a useful characterization of *r*-statistically convergent sequences.

**Theorem 2.2.** Let $r \in [0,1]$. Then, $a = r$-*stat*-lim $l$ if and only if there exists an increasing index sequence $K = \{k_n\,;\,k_n \in N, n = 1, 2, 3, …\}$ of the natural numbers such that $d(K) = 1 - r$ and $a = \lim l_K$ where $l_K = \{a_i\,;\,i \in K\}$.

<u>Proof</u>. *Necessity*. Let us consider a sequence $l = \{a_i\,;\,i = 1, 2, 3, …\}$ of real numbers, for which $a = r$-*stat*-lim $l$. Then for any natural number $n$, we can build the set $L_{1/2n}(a) = \{i \in N;\,|a_i - a| \geq 1/2n\}$. As $a = r$-*stat*-lim $l$, we have $0 = d_r(L_\varepsilon(a))$. It means that there is a number $p$ such that $(1/p)\,|\,L_{1/2n}(a)\,| < r + 1/2n$. As the number $p$ depends on $n$, we denote it by $p(n)$.

Using the sequence of numbers $p(n)$, we build a sequence $\{l_0, l_1, …, l_n, …\}$ of subsequences of the sequence $l$. Here $l_0 = l$, $l_1 = l_0 \setminus \{a_i\,;\,i \in L_{1/2n}(a)\}$, …, $l_n = l_{n-1} \setminus \{a_i\,;\,i \in L_{1/2n}(a)\}$, …, $n = 1, 2, 3, …$ Let us put $h = \bigcap_{n=1}^{\infty} l_n$.

By our construction, the ratio $(1/p(n))\,|L_{1/2n}(a)|$ converges to $r$ when $n$ tends to infinity. Let us denote by $I(h)$ indices of all those element from $l$ that belong to $h$. Then $d(I(h)) = 1 - r$ and $a = stat\text{-}\lim h$.

By Theorem 3.2 from (see Burgin and Duman [4]) $a = stat\text{-}r\text{-}\lim l$ if and only if there exists an increasing index sequence $K = \{k_n\,;\,k_n \in N,\,n = 1,2,3,\ldots\}$ of the natural numbers such that $d(K) = 1$ and $a = r\text{-}\lim l_K$ where $l_K = \{a_i\,;\,i \in K\}$. Thus, in the case of $r = 0$, there is a subsequence $k = \{a_i\,;\,i \in J\}$ of $h$ such that $J \subseteq I(h)$ and $d_{I(h)}(J) = \lim_{n\to\infty} |J_n|/n = 1$, where $J_n$ are those numbers of the elements from $h$ that belong to $k$ and to the first $n$ elements from $h$. By our construction, we have $a = \lim k$.

Necessity is proved.

*Sufficiency.* Let us consider a sequence $l = \{a_i\,;\,i = 1, 2, 3,\ldots\}$ of real numbers and its subsequence $k = \{a_j\,;\,j \in J \subseteq N\}$ of $h$ such that $d(J) = 1 - r$ and $a = \lim k$. Then $d(N \setminus J) = 1 - (1 - r) = r$. This means that 0 is an $r$-statistical limit of $l$.

Theorem is proved.

Let $l = \{a_i \in R;\,i = 1, 2, 3,\ldots\}$ and $h = \{b_i \in R;\,i = 1, 2, 3,\ldots\}$. Then their sum $l + h$ is equal to the sequence $\{a_i + b_i;\,i = 1, 2, 3,\ldots\}$ and their difference $l - h$ is equal to the sequence $\{a_i - b_i;\,i = 1, 2, 3,\ldots\}$. Lemma 2.6 allows us to prove the following result.

**Theorem 2.3.** Let $a = r\text{-}stat\text{-}\lim l$ and $b = q\text{-}stat\text{-}\lim h$. Then we have

**(a)** $a + b = (r+q)\text{-}stat\text{-}\lim(l+h)$;

**(b)** $a - b = (r+q)\text{-}stat\text{-}\lim(l - h)$;

**(c)** $ka = (r\,|k|)\text{-}stat\text{-}\lim(kl)$ for any $k \in R$ where $kl = \{ka_i\,;\,i = 1,2,3,\ldots\}$.

<u>Proof</u>. **(a)** Let $l = \{a_i\,;\,i = 1,2,3,\ldots\}$ and $h = \{b_i\,;\,i = 1,2,3,\ldots\}$. Then, by hypotheses, we get, for a given $\varepsilon > 0$, that $0 = r\text{-}\lim_{n\to\infty}(1/n)\,|\{i \in N;\,i \le n\,;\,|a_i - a| \ge \varepsilon/2\}|$ and $0 = q\text{-}\lim_{n\to\infty}(1/n)\,|\{i \in N;\,i \le n\,;\,|b_i - b| \ge \varepsilon/2\}|$. Set

$$u_n = |\{i \in N;\,i \le n\,;\,|a_i - a| \ge \varepsilon/2\}| \text{ and } v_n = |\{i \in N;\,i \le n\,;\,|b_i - b| \ge \varepsilon/2\}|.$$

Now observe that

$$|\{i \in N;\,i \le n\,;\,|(a_i + b_i) - (a + b)| \ge \varepsilon\}| \le u_n + v_n,$$

which implies

$$(1/n) \, |\{i \in N; i \leq n \, ; \, |(a_i + b_i) - (a + b)| \geq \varepsilon\}| \leq (1/n) \, u_n + (1/n) \, v_n.$$

Thus, for any $\alpha > 0$ and for sufficiently large $n$, we may write that

$$(1/n) \, |\{i \in N; i \leq n \, ; \, |(a_i + b_i) - (a + b)| \geq \varepsilon\}| \leq r + \alpha + q + \alpha = (r+q) + 2\alpha.$$

Therefore, we obtain that

$$a + b = (r+q)\text{-}stat\text{-}\lim(l+h).$$

Proofs of the statements **(b)** and **(c)** are similar.

**Corollary 2.4** [16]. If $b = stat\text{-}\lim l$ and $c = stat\text{-}\lim h$, then:

**(a)** $a + b = stat\text{-}\lim (l+h)$;

**(b)** $a - b = stat\text{-}\lim (l - h)$;

**(c)** $ka = stat\text{-}\lim (kl)$ for any $k \in R$.

Now we can obtain the following "squeeze rule" for $r$-statistically convergent sequences. Let $l = \{a_i \in R; i = 1, 2, 3, \ldots\}$, $h = \{b_i \in R; i = 1, 2, 3, \ldots\}$ and $k = \{c_i \in R; i = 1, 2, 3, \ldots\}$.

**Theorem 2.4.** If $a_i \leq b_i \leq c_i$ ( $i = 1, 2, 3, \ldots$ ), and if $a = r\text{-}stat\text{-}\lim l = r\text{-}stat\text{-}\lim k$, then we have $a = r\text{-}stat\text{-}\lim h$.

Proof. For a given $\varepsilon > 0$, let

$$u_n = |\{i \in N; i \leq n \, ; \, -\varepsilon < a_i - a\}|,$$

$$v_n = |\{i \in N; i \leq n \, ; \, |b_i - a| < \varepsilon\}|,$$

$$y_n = |\{i \in N; i \leq n \, ; \, c_i - a < \varepsilon\}|.$$

Since $a_i \leq b_i \leq c_i$ ($i = 1, 2, 3, \ldots$ ), we easily see that $u_n/n \leq v_n/n \leq y_n/n$ for each $n \in N$. From the initial conditions and Lemma 2.6 (b), we know that $1 = r\text{-}\lim (u_n/n) = r\text{-}\lim (y_n/n)$. Thus, for any $\alpha > 0$ and sufficiently large $n$, we may write that

$$-r - \alpha < (u_n/n) - 1 \leq (v_n/n) - 1 \leq (y_n/n) - 1 < r + \alpha,$$

which guarantees that $|(v_n/n) - 1| < r + \alpha$ for any $\alpha > 0$ and sufficiently large $n$. Consequently, we get $1 = r\text{-}\lim (v_n/n)$, which yields, by Lemma 2.6, that $a = r\text{-}stat\text{-}\lim h$.

The proof is completed.

### 3. Fuzzy convergence of statistical characteristics

To each sequence $l = \{a_i\, ; i = 1, 2, 3, \ldots\}$ of real numbers, it is possible to correspond a new sequence $\mu(l) = \{\mu_n = (1/n) \sum_{i=1}^{n} a_i\, ; n = 1,2,3,\ldots\}$ of its partial averages (means). Here a partial average of $l$ is equal to $\mu_n = (1/n) \sum_{i=1}^{n} a_i$. Partial averages are also called Cesáro means. They are used to determine Cesáro summability in the theory of divergent series. Cesáro summability is often applied to Fourier series as it is more powerful than traditional summability.

Sequences of partial averages/means play an important role in the theory of ergodic systems [1]. Indeed, the definition of an ergodic system is based on the concept of the "time average" of the values of some appropriate function $g$ arguments for which are dynamic transformations $T$ of a point $x$ from the manifold of the dynamical system. This average is given by the formula

$$\hat{g}(x) = \lim (1/n) \sum_{k=1}^{n-1} g(T^k x).$$

In other words, the dynamic average is the limit of the partial averages/means of the sequence $\{T^k x\, ; k = 1,2,3,\ldots\}$.

Properties of the average $\hat{g}(x)$ when the parameter $k$ is a discrete time are described in the famous Birkhoff-Khinchin theorem, which is one of the most important results in ergodic theory [17].

Let $l = \{a_i\, ; i = 1, 2, 3, \ldots\}$ be a bounded sequence, i.e., there is a number $m$ such that $|a_i| < m$ for all $i \in N$. This condition is usually true for all sequences generated by measurements or computations, i.e., for all sequences of data that come from real life.

**Theorem 3.1.** If $a = r\text{-}stat\text{-}\lim l$, then $a = u\text{-}\lim \mu(l)$, where $u = (m + |a|)r$.

Proof. Since $a = r\text{-}stat\text{-}\lim l$, for every $\varepsilon > 0$, we have

$$0 = r\text{-}\lim_{n\to\infty} (1/n) |\{i \leq n, i \in N; |a_i - a| \geq \varepsilon\}| \qquad (3.1)$$

If $|a_i| < m$ for all $i \in N$, then there is a number $k$ such that $|a_i - a| < k$ for all $i \in N$. Namely, $|a_i - a| \leq |a_i| + |a| \leq m + |a| = k$. Taking the set $L_{n,\varepsilon}(a) = \{i \leq n, i \in N; |a_i - a| \geq \varepsilon\}$, denoting $|L_{n,\varepsilon}(a)|$ by $u_n$, and using the hypothesis $|a_i| < m$ for all $i \in N$, we have the following system of inequalities:

$$|\mu_n - a| = |(1/n) \sum_{i=1}^{n} a_i - a|$$
$$\leq (1/n) \sum_{i=1}^{n} |a_i - a|$$
$$\leq (1/n) \{ku_n + (n - u_n)\varepsilon\}$$
$$\leq (1/n) (ku_n + n\varepsilon)$$
$$= \varepsilon + k (u_n/n).$$

From (3.1), we get the inequality $|\mu_n - a| < \varepsilon(1 + k) + kr$ for sufficiently large $n$. Since $m + |a| = k$, we have $a = u\text{-}\lim \mu(l)$, where $u = (m + |a|)r$.

Theorem is proved.

Taking $r = 0$, we have the following result.

Let $l = \{a_i ; i = 1,2,3,\ldots\}$ be a bounded sequence such that $|a_i| < m$ for all $i \in N$.

**Corollary 3.1.** If $a = stat\text{-}\lim l$, then $a = \lim \mu(l)$.

**Remark 3.1.** However, (fuzzy) convergence of the partial averages/means of a sequence does not imply (fuzzy) statistical convergence of this sequence.

**Remark 3.2.** For unbounded sequences, the result of Theorem 3.1 can be invalid. Thus, the condition of boundedness is essential in this theorem.

Taking the sequence $l = \{a_i ; i = 1, 2, 3,\ldots\}$ of real numbers, it is possible to construct not only the sequence $\mu(l) = \{ \mu_n = (1/n) \sum_{i=1}^{n} a_i ; n = 1, 2, 3,\ldots\}$ of its partial averages (means) but also the sequences $\sigma(l) = \{\sigma_n = ((1/n) \sum_{i=1}^{n} (a_i - \mu_n)^2)^{1/2} ; n = 1,2,3,...\}$ of its partial standard deviations $\sigma_n$ and $\sigma^2(l) = \{\sigma_n^2 = (1/n) \sum_{i=1}^{n} (a_i - \mu_n)^2 ; n = 1,2,3,...\}$ of its partial variances $\sigma_n^2$.

Let us find how fuzzy statistical convergence of a sequence is related to fuzzy statistical convergence of this sequence of its partial standard deviations.

**Theorem 3.2.** If $a = r$-stat-lim $l$ and $|a_i| < m$ for all $i = 1, 2, 3, \ldots$, then $0 = [p(2r + u)]^{1/2}$-lim $\sigma(l)$ where $p = \max\{m^2 + |a|^2, m + |a|\}$ and $u = (m + |a|)r$.

<u>Proof.</u> At first, we show that $\lim \sigma^2(l) = 0$. By the definition, $\sigma_n^2 = (1/n) \sum_{i=1}^{n} (a_i - \mu_n)^2 = (1/n) \sum_{i=1}^{n} (a_i)^2 - \mu_n^2$. Thus, $\lim \sigma^2(l) = \lim_{n\to\infty} (1/n) \sum_{i=1}^{n} (a_i)^2 - \lim_{n\to\infty} \mu_n^2$. Since $|a_i| < m$ for all $i \in N$, there is a number $p$ such that $|a_i^2 - a^2| < p$ for all $i \in N$. Namely, $|a_i^2 - a^2| \leq |a_i|^2 + |a|^2 < m^2 + |a|^2 < \max\{m^2 + |a|^2, m + |a|\} = p$. Taking the set $L_{n,\varepsilon}(a) = \{i \in N; i \leq n \text{ and } |a_i - a| \geq \varepsilon\}$, denoting $|L_{n,\varepsilon}(a)|$ by $u_n$, and using the hypothesis $|a_i| < m$ for all $i \in N$, we have the following system of inequalities:

$$|\sigma^2_n| = |(1/n) \sum_{i=1}^{n} (a_i)^2 - \mu_n^2|$$
$$= |(1/n) \sum_{i=1}^{n} (a_i^2 - a^2) - (\mu_n^2 - a^2)|$$
$$\leq (1/n) \sum_{i=1}^{n} |a_i^2 - a^2| + |\mu_n^2 - a^2|$$
$$< (p/n) \sum_{i=1}^{n} |a_i - a| + |\mu_n - a| |\mu_n + a|$$
$$< (p/n) (u_n + (n - u_n)(r + \varepsilon)) + |\mu_n - a| (|\mu_n| + |a|)$$
$$\leq (p/n) (u_n + n(r + \varepsilon)) + |\mu_n - a| ((1/n)\sum_{i=1}^{n} |a_i| + |a|)$$
$$< p(u_n/n) + p(r + \varepsilon) + p|\mu_n - a|.$$

Now by hypothesis and Theorem 3.1, we have $a = u$-lim $\mu(l)$, where $u = (m + |a|)r$. Also, for every $\varepsilon > 0$ and sufficiently large $n$, we may write that

$$|\sigma^2_n| < p\varepsilon + pr + p(r + \varepsilon) + p(u + \varepsilon) = p(2r + u) + 3p\varepsilon \qquad (3.2)$$

As $(x + y)^{1/2} \leq x^{1/2} + y^{1/2}$ for any $x, y > 0$, the inequality (3.2) implies the inequality

$$|\sigma_n| \leq [p(2r + u)]^{1/2} + (3p\varepsilon)^{1/2},$$

which yields that $0 = [p(2r + u)]^{1/2}$-lim $\sigma(l)$.

The proof is completed.

Taking *r* = 0, we have the following result for bounded sequences.

**Corollary 3.2.** If *a* = *stat*-lim *l*, then 0 = lim σ(*l*).

### 4. Conclusion

We have developed the concept of fuzzy statistical convergence and studied its properties. In particular, it is demonstrated that fuzzy statistical convergence of a sequence implies fuzzy convergence of partial averages (means) of this sequence (Theorem 3.1) and fuzzy convergence of partial standard deviations of this sequence (Theorem 3.2).

Some relations between fuzzy statistical convergence and fuzzy summability are obtained for Cesáro summability. However, it would be interesting to study relations between fuzzy statistical convergence and fuzzy summability in more detail.

Results obtained in this paper also open an approach to the development of fuzzy ergodic theory. In particular, it would be interesting to extend the Birkhoff-Khinchin theorem for fuzzy measures and fuzzy limits.